\documentclass[12pt]{amsart}

\usepackage{latexsym,amsfonts,amssymb,epsfig,verbatim}
\usepackage{amsmath, latexsym, graphics, textcomp}
\usepackage{amsthm,diagbox}
\usepackage{pinlabel}

\usepackage[all,2cell,dvips]{xy}

\newcommand{\nc}{\newcommand}
\nc{\dmo}{\DeclareMathOperator}
\nc{\nt}{\newtheorem}

\nt{theorem}{Theorem}[section]
\nt{lemma}[theorem]{Lemma}
\nt{prop}[theorem]{Proposition}
\nt{question}[theorem]{Question}
\nt{cor}[theorem]{Corollary}
\nt{fact}[theorem]{Fact}
\nt{problem}[theorem]{Problem}
\nt{conjecture}[theorem]{Conjecture}
\nt{meta}[theorem]{Metaconjecture}
\nt{genmeta}[theorem]{Generalized Metaconjecture}

\nc{\Z}{\mathbb{Z}}
\nc{\R}{\mathbb{R}}
\nc{\Q}{\mathbb{Q}}
\nc{\C}{\mathcal{C}}
\nc{\complex}{\mathbb{C}}
\nc{\I}{\mathcal{I}}
\nc{\BI}{\mathcal{BI}}
\nc{\K}{\mathcal{K}}
\nc{\N}{\mathcal{N}}

\dmo{\GL}{GL}
\dmo{\SL}{SL}
\dmo{\Sp}{Sp}
\dmo{\Mod}{Mod}
\dmo{\B}{B}
\dmo{\PB}{PB}
\dmo{\Aut}{Aut}
\dmo{\Out}{Out}
\dmo{\UT}{UT}

\newcommand{\BigFreeProd}{\mathop{\mbox{\Huge{$\ast$}}}}

\nc{\margin}[1]{\marginpar{\tiny #1}}

\nc{\p}[1]{\smallskip\noindent{{\bf #1}}}

\include{diagram}

\begin{document}


\title[Problems, Questions, and Conjectures]{Problems, Questions, and Conjectures about Mapping Class Groups}

\author{Dan Margalit}

\address{School of Mathematics, 686 Cherry St, Atlanta, GA 30306}

\thanks{This material is based upon work supported by the National Science Foundation under Grant No. DMS - 1057874}



\begin{abstract}
We discuss a number of open problems about mapping class groups of surfaces.  In particular, we discuss problems related to linearity, congruence subgroups, cohomology, pseudo-Anosov stretch factors, Torelli subgroups, and normal subgroups.
\end{abstract}

\maketitle

Beginning with the work of Max Dehn a century ago, the subject of mapping class groups has become a central topic in mathematics.  It enjoys deep and varied connections to many other subjects, such as low-dimensional topology, geometric group theory, dynamics, Teichm\"uller theory, algebraic geometry, and number theory.   The number of papers on mapping class groups recorded on MathSciNet in the last six decades has grown from 205 to 386 to 525 to 791 to 1,121 to 1,390.  The subject seems to enjoy an endless supply of beautiful ideas, pictures, and theorems.  

At the 2017 Georgia International Topology Conference, the author gave a lecture called ``Problems and progress on mapping class groups.''  This paper is a summary of parts of that lecture.  What follows is not a comprehensive list in any way, even among the topics it attempts to address.  Rather it gives a mix of problems---from the famous and notoriously difficult to the eminently doable.  There is little attempt to give background; the reader may find that in the book by Farb and the author \cite{primer} and in the other references therein.

There are other problem lists on mapping class groups.  In fact Farb has edited an entire book of problem lists on mapping class groups \cite{problembook}.  Ivanov's problem list \cite{15}, which also appears in Farb's book, has been particularly influential on the author of this article.  That problem list is an updated version of a problem list written in conjunction with the 1993 Georgia International Topology Conference, twenty-four years earlier.  Some of Ivanov's problems also appeared in Kirby's problem list, compiled around the same time \cite{kirby}.  Joan Birman's classic book contains several problems on mapping class groups, some of which are still open \cite{jsb}.

For a surface $S$, its mapping class group $\Mod(S)$ is the group of homotopy classes of orientation-preserving homeomorphisms of $S$.  If $S$ has boundary and/or marked points, then the homeomorphisms (and homotopies) are required to fix the marked points as a set and are required to fix the boundary of $S$ pointwise.  We denote the closed orientable surface of genus $g$ by $S_g$.

\subsection*{Acknowledgments} The author would like to thank Tara Brendle, Lei Chen, Nathan Dunfield, Benson Farb, Autumn Kent, Eiko Kin, Kevin Kordek, Justin Lanier, Johanna Mangahas, Andy Putman, Nick Salter, Shane Scott, Bal\'azs Strenner, and the anonymous referee.

\section{Linearity}

We begin with one of the most basic and famous questions about mapping class groups, which appears, for example, in Birman's book \cite[Appendix, Problem 30]{jsb}.  In the statement, we say that a group is linear if it admits a faithful linear representation into some $\GL_n(k)$ where $k$ is a field.  

\begin{question}
\label{linear}
Is $\Mod(S_g)$ linear?
\end{question}

Dehn proved that $\Mod(T^2)$, the mapping class group of the torus, is isomorphic to $\SL_2(\Z) \subseteq \GL_2(\complex)$.  In the case $g=2$ the Birman--Hilden theory \cite{bht} gives a short exact sequence:
\[
1 \to \langle \iota \rangle \to \Mod(S_2) \to \Mod(S_{0,6}) \to 1
\]
where $\iota$ is the hyperelliptic involution of $S_2$.  The group $\Mod(S_{0,6})$ is closely related to the braid group on 5 strands.  As such, Bigelow--Budney and Korkmaz were able to prove that $\Mod(S_2)$ is linear, using the theorem of Krammer and Bigelow that braid groups are linear \cite{BigelowLinear,KrammerLinear}.  For $g \geq 3$,  Question~\ref{linear} is wide open.

\subsection*{Linearity of the braid group} One might be tempted to think that $\Mod(S_g)$ is not linear, because if it were then we would already have stumbled across the representation.  On the other hand, we should draw inspiration from the case of the braid group.  The Burau representation of the braid group was introduced in 1935 \cite{Burau}.  It was shown to be unfaithful for $n \geq 6$ by Long--Paton and Moody \cite{LongMoody,moody} in 1993 and later it was shown to be unfaithful for $n \geq 5$ by Bigelow \cite{Bigelow5}.  (The case $n=4$ is open; see Question~\ref{b4} below) The Lawrence--Krammer representation was introduced in 1990 by Lawrence \cite{Lawrence} and was proved to be faithful by Bigelow and Krammer around the turn of the century.  See the surveys by Birman and Brendle \cite{braids} and by Turaev \cite{turaev}.

Whereas the Burau representation considers the action of the braid group on the fundamental group of the punctured disk, the Lawrence--Krammer representation considers the action on the fundamental group of the space of configurations of two points in the punctured disk.  The fundamental group of a space is the fundamental group of the space of configurations of a single point, and so in that sense the Lawrence--Krammer representation is a mild generalization.  This gives hope that some mild generalization of a known representation of $\Mod(S_g)$ could also be faithful.

\subsection*{Lawrence--Krammer representations for mapping class groups} Let us define the $k$th Lawrence representation of $\Mod(S_g)$ to be the action of $\Mod(S_g)$ on the homology of the universal abelian cover of the configuration space of $k$ distinct, ordered points in $S_g$, considered as a module over $H_1(S_g;\Z)$.  (This is not a perfect analog of the Lawrence--Krammer representation because there is no $\Mod(S_g)$-equivariant surjective map $H_1(S_g;\Z) \to \Z$.)

\begin{question}
\label{law}
Is the $k$th Lawrence representation of $\Mod(S_g)$ faithful for any $k \geq 1$?
\end{question}

Some positive evidence\footnote{According the the anonymous referee, Moriyama's result should be considered as negative evidence!} in the direction of Question~\ref{law} is given by the work of Moriyama \cite{moriyama}, who showed that the $k$th term of the Johnson filtration for the once-punctured surface $S_{g,1}$ is given by the kernel of the action of $\Mod(S_{g,1})$ on the compactly supported cohomology of the configuration space of $k$ distinct points in $S_{g,1}$. 

\subsection*{Kontsevich's approach} Kontsevich has proposed a family of linear representations of the mapping class group, described in the online problem list for the Center for the Topology and Quantization of Moduli Spaces \cite{CTQM} (there is also a discussion on MathOverflow \cite{MO}).  Kontsevich has conjectured that this approach gives a faithful linear representation of the mapping class group.  Here is his description of the main idea, taken from the CTQM web site and mildly edited for clarity.

\smallskip

\emph{According to Thurston, a generic element $f$ in $\Mod(S_g)$ is pseudo-Anosov and has a canonical representative that  preserves two transverse laminations and multiplies the transverse measures by $\lambda$ and $1/\lambda$, where $\lambda>1$ is an algebraic integer.  The number $\lambda$ appears as the largest eigenvalue of the action of $f$ on the first cohomology of the ramified double covering of $S_g$ on which the quadratic differential associated to the laminations for $f$ becomes an abelian differential. Moreover, one may assume that the laminations have $4g-4$ triple singular points. I would like to construct a representation in which $f$ acts non-trivially, because there will be some eigenvalue of the corresponding matrix equal to a positive power of $\lambda$.}

\emph{A representation of $\Mod(S_g)$ is the same as a local system on moduli space $\mathcal{M}_g$. Here is the description of it. The fiber of it at a complex curve $C$ is a kind of middle (or total?) cohomology of the configuration space $\textrm{Conf}_{4g-4}(C)=C^{4g-4}\setminus\textrm{big diagonal}$ (the configuration space of $4g-4$ distinct points $(x_1,...,x_{4g-4})$ on $C$), with coefficients in the flat bundle, the first cohomology group of the double cover of $C$, doubly ramified at $(x_1,...,x_{4g-4})$.}

\emph{Why I believe that this should work: Thurston's representative gives a map from $\textrm{Conf}_{4g-4}(C)$ (together with the local system) to itself, with a preferred marked point such that the action on the fiber has the largest eigenvalue equal to $\lambda$.
Can one use Lefschetz fixed point formula, or maybe some dynamical reasoning, to see that this eigenvalue appear in the total cohomology of $\textrm{Conf}_{4g-4}(C)$ with local coefficients as well?}

\smallskip

There are various challenges to making Kontsevich's idea work.  For instance, it is not clear why one can assume that the laminations only have singular points of degree 3; it is true that any foliation is Whitehead equivalent to such a foliation, but then there is no pseduo-Anosov homeomorphism preserving the modified foliation.  Even assuming that all pseudo-Anosov foliations have only singular points of degree 3, it is still not clear how to construct the local system; it seems that this would require a deformation retraction of $\textrm{Conf}_{4g-4}(C)$ to some base point. 

\begin{question}
Can Kontsevich's idea be fleshed out to give a faithful linear representation of the mapping class group, or not?
\end{question}

In unpublished work, Dunfield has performed some computer calculations suggesting that the Kontsevich approach will not result in a faithful representation, at least for the braid groups \cite{dunfield}.

\subsection*{No poison} We end this section by mentioning one other result in the direction of Question~\ref{linear}.  Let $F_n$ denote the free group of rank $n$.  Formanek--Procesi \cite{fp} showed that the automorphism group $\Aut(F_n)$ and the outer automorphism group $\Out(F_n)$ both fail to be linear for $n \geq 4$.  They accomplish this  by producing certain ``poison subgroups'' that are patently non-linear (the $n=3$ cases are open).  Brendle and Hamidi-Tehrani showed that there are no analogous poison subgroups in the mapping class group \cite{poison}.


\section{The congruence subgroup problem}
\label{cong}

For a given group, the congruence subgroup problem asks if every finite-index subgroup of a given group contains a congruence subgroup.  The precise meaning of this problem depends on an appropriate choice of definition for a congruence subgroup.

The \emph{principal congruence subgroup of level $m$} in the special linear group $\SL_n(\Z)$ is the kernel of the reduction homomorphism $\SL_n(\Z) \to \SL_n(\Z/m)$.  It is a theorem of Bass--Milnor--Serre \cite{bms} that $\SL_n(\Z)$ has the congruence subgroup property.  As above, this means that every subgroup of finite index in $\SL_n(\Z)$ contains a principal congruence subgroup.

The most na\"ive analogue of the principal congruence subgroups for $\Mod(S_g)$ are the level $m$ congruence subgroups $\Mod(S_g)[m]$, defined to be the kernels of the compositions
\[
\Mod(S_g) \to \Sp_{2g}(\Z) \to \Sp_{2g}(\Z/m),
\]
where $\Sp_{2g}(\Z)$ and $\Sp_{2g}(\Z/m)$ are the symplectic groups over $\Z$ and $\Z/m$.  This composition gives exactly the action of $\Mod(S_g)$ on the homology group $H_1(S_g;\Z/m)$.  The Torelli group $\I(S_g)$ is the kernel of the map $\Mod(S_g) \to \Sp_{2g}(\Z)$, and as such it is the intersection of all of the level $m$ congruence subgroups for $m \geq 2$.  

There are finite-index subgroups of $\Mod(S_g)$ not containing any congruence subgroup $\Mod(S_g)[m]$.  A fortiori there are subgroups of finite index that do not contain the Torelli group.  Here is how we can construct such a subgroup.  Let $H$ denote $H_1(S_g;\Z)$.  Morita \cite[Theorem 6.2]{moritaext} defined a homomorphism
\[
\tilde \tau : \Mod(S_g) \to \left(\wedge^3 H / H\right) \rtimes \Sp_{2g}(\Z) 
\]
that extends the Johnson homomorphism $\tau : \I(S_g) \to \wedge^3 H / H$, in the sense that the restriction of $\tilde \tau$ to $\I(S_g)$ is equal to the post-composition of $\tau$ with the inclusion $\left(\wedge^3 H / H\right) \to \left(\wedge^3 H / H\right) \rtimes \Sp_{2g}(\Z)$.  Consider the subgroup $2\wedge^3 H / 2H$ where the inclusion $2H \to 2\wedge^3 H$ is the restriction of the usual inclusion $H \to \wedge^3 H$.  The action of $\Sp_{2g}(\Z)$ on $\wedge^3 H$ preserves both $2\wedge^3 H$ and $2H$, and so we obtain a subgroup $\left(2\wedge^3 H / 2H\right) \rtimes \Sp_{2g}(\Z)$ inside $\left(\wedge^3 H / H\right) \rtimes \Sp_{2g}(\Z)$.  The preimage of this subgroup under $\tilde \tau$ is the desired subgroup.  

We may vary this construction in the following way: if $\Gamma$ is a subgroup of finite index in $\Sp_{2g}(\Z)$ and $\Delta$ is a subgroup of finite index in $\wedge^3 H / H$ invariant under $\Gamma$, then $\Delta \rtimes \Gamma$ is a finite-index subgroup of $\Mod(S_g)$.  In what follows we will refer to these as Morita subgroups.

Other examples of finite-index subgroups of $\Mod(S_g)$ that do not contain the Torelli group $\I(S_g)$ are given by Cooper \cite{cooper}, who constructs two different characteristic $p$ analogues of the Johnson filtration. 

The upshot is that if we want the mapping class group to have the congruence subgroup property then we shall require a richer class of subgroups of finite index.

To this end, consider a characteristic subgroup $\Gamma$ of finite index in $\pi_1(S_g)$.  There is a homomorphism $\Mod(S_g) \to \Out (\pi_1(S_g)/\Gamma)$.  Since $\Gamma$ has finite index the group $\pi_1(S_g)/\Gamma$, hence $\Out (\pi_1(S_g)/\Gamma)$, is finite.  The kernel of the homomorphism is thus a subgroup of finite index.  

For instance if $\Gamma$ is the kernel of the composition $\pi_1(S_g) \to H_1(S_g;\Z) \to H_1(S_g;\Z/m)$ then we recover the congruence subgroup $\Mod(S_g)[m]$.  Following Ivanov, we refer to any subgroup of $\Mod(S_g)$ constructed in this manner as a congruence subgroup.

\begin{question}
\label{cong}
Is it true that every finite-index subgroup in $\Mod(S_g)$ contains a congruence subgroup?
\end{question}

This is the first problem on Ivanov's list \cite{15}, and in fact he states it as a conjecture.  More modestly, we may ask the following.

\begin{question}
Do the Morita subgroups of $\Mod(S_g)$ constructed above contain congruence subgroups?
\end{question}

It would even be helpful to clarify the relationship between the Morita subgroups and the Cooper subgroups.

We can reinterpret the construction of congruence subgroups in terms of covering spaces.  Fix a characteristic subgroup $\Gamma$ of finite index in $\pi_1(S_g)$ and consider the associated characteristic finite-sheeted cover of $S_g$.  The entire mapping class group $\Mod(S_g)$ is liftable in the sense of Birman--Hilden \cite{bh,bht}, that is, each element lifts to a mapping class of the cover.  As such, $\Mod(S_g)$ acts by outer automorphisms on the deck group of the cover, and the kernel of this action is precisely the congruence subgroup associated to $\Gamma$.

This idea gives another source for subgroups of finite index.  For a subgroup $\Gamma$ of finite index in $\pi_1(S_g)$, we may consider the associated finite-sheeted cover of $S_g$.  The resulting liftable mapping class group is again a subgroup of $\Mod(S_g)$ of finite index.  It is a straightforward exercise to show that these liftable mapping class groups all contain congruence subgroups as defined by Ivanov.

There are some surfaces for which Question~\ref{cong} has been resolved.  Diaz--Donagi--Harbater \cite{DDH} answered it in the affirmative for a sphere with four punctures (although their proof applies to spheres with any number of punctures).  Asada \cite{asada} answered the question in the affirmative for the sphere with any number of punctures and for the torus with any positive number of punctures (see also \cite{BER,EM}).  Thurston gave a more elementary proof for punctured spheres, based on Asada's work; see the paper by McReynolds \cite{ben} for a presentation of Thurston's proof as well as some historical discussion.  More recently, a new proof in the case of the punctured torus was given by Kent \cite{aek}.  Question~\ref{cong} was also resolved for surfaces of genus two by Boggi \cite{boggi,boggi2}.  (Boggi \cite{boggi3} had claimed a general solution to Question~\ref{cong} but a gap was found by Abramovich, Kent, and Wieland in his Theorem 5.4.)

As mentioned by Ivanov \cite{15}, Voevodsky observed that a positive answer to Question~\ref{cong} would give a proof of a conjecture of Grothendieck that a smooth algebraic curve defined over $\Q$ is determined by its algebraic fundamental group---which is isomorphic to the profinite completion of $\pi_1(S_g)$---together with the natural action of the absolute Galois group $\textrm{Gal}(\bar \Q /\Q)$ on the latter (the conjecture was proved by Mochizuki \cite{mochi1,mochi2}).


\section{The integral Burau representation}

Denote by $\B_n$ the braid group on $n$ strands.  The Burau representation
\[
\B_n \to \GL_n(\Z[t,t^{-1}])
\]
is the most well-studied linear representation of $\B_n$.  It is easily seen to be faithful for $n=3$.  As discussed above, it is now known to be unfaithful for $n \geq 5$.  This begs the following question, discussed, for example, in Birman's book \cite{jsb}.

\begin{question}
\label{b4}
Is the Burau representation faithful for $n=4$?
\end{question}

Bigelow suggests that an element of the kernel of the Burau representation for $n=4$ would give rise to a candidate knot for showing that the Jones polynomial does not detect the unknot \cite{bigjones}.  Some experimental results on Question~\ref{b4} were recently found by Fullarton and Shadrach \cite{fs}.  Birman showed that Question~\ref{b4} is equivalent to the question of whether the images of $\sigma_1\sigma_3^{-1}$ and $\sigma_2\sigma_3\sigma_1^{-1}\sigma_2^{-1}$ generate a free group of rank 2 \cite[Theorem 3.19]{jsb}.

While the Burau representation is not faithful in general, we can still use it to effectively probe the braid group.  For instance the kernel of the Burau representation gives an interesting normal subgroup of the braid group.  Church and Farb showed that this group is not finitely generated \cite{cf}.  Still, we have the following problem, implicit in Morita's survey \cite[Problem 6.24]{prospect}.

\begin{problem}
Describe the kernel of the Burau representation.  Find a natural generating set.
\end{problem}

In Morita's survey, he mentions that it is still open whether the Gassner representation of the pure braid group is faithful or not; see Gassner's paper \cite{gassner} for the definition.  Bachmuth \cite{bachmuth} had claimed that it was faithful, but his argument was refuted by Abramenko--M\"uller \cite{am}.  We refer the reader to Birman's extensive review of Bachmuth's paper on MathSciNet.

\subsection*{Braid Torelli groups} We may also consider specializations of the Burau representation obtained by setting $t$ to be some fixed complex number.  There are two values of $t$ that give specializations to $\GL_n(\Z)$, namely $t=1$ and $t=-1$.  The specialization at $t=1$ is nothing other than the standard representation to the symmetric group.  The specialization at $t=-1$, called the integral Burau representation, is more interesting.  (Specializations at other roots of unity have also been studied; see, e.g. \cite{bb,fk,gg,venky}.)

For $n=2g+1$ there is a map $\B_{2g+1} \to \Mod(S_g^1)$ obtained by lifting homeomorphisms through the branched double cover $S_g^1 \to D_{2g+1}$ (here $S_g^1$ is the compact surface of genus $g$ with one boundary component and $D_{2g+1}$ is the disk with $2g+1$ marked points).  The integral Burau representation is the representation obtained by post-composing this map with the standard symplectic representation $\Mod(S_g^1) \to \Sp_{2g}(\Z)$.  For $n$ even the integral Burau representation has a similar description, with $S_g^1$ replaced by $S_g^2$;  see \cite{BMP}.  

The kernel of the integral Burau representation is called the braid Torelli group $\BI_n$.  Hain has identified this group with the fundamental group of the branch locus for the period mapping on Torelli space \cite{hainfinite}, and so this group is of interest in algebraic geometry.  With Brendle and Putman, the author of this article showed that $\BI_n$ is generated by squares of Dehn twists about curves in $D_{2g+1}$ that surround an odd number of marked points \cite{BMP}.  

Smythe \cite{smythe} showed that $\BI_4$ is isomorphic to $F_\infty$, the free group on countably many generators.  It follows from a theorem of Mess \cite{Mess} that $\BI_5 \cong F_\infty \times \Z$ and then from a theorem of Brendle and the author \cite{pp} that $\BI_6 \cong F_\infty \ltimes F_\infty$.  The following problem is suggested by Hain \cite[Problem 7]{hainfinite}, although it was clearly of interest much earlier.

\begin{question}
Is $\BI_n$ finitely generated for any $n \geq 7$?   Is the abelianization finitely generated?
\end{question}

It is further of interest to find (possibly infinite) presentations of the groups $\BI_n$.  

\subsection*{Level $m$ congruence subgroups} We obtain finite-index analogs of $\BI_n$ by considering the mod $m$ reductions of the integral Burau representation.  The kernels of these representations are called the level $m$ congruence subgroups of the braid group, and are denoted $\B_n[m]$.  Arnol'd \cite{arnold} proved that $\B_n[2]$ is nothing other than the pure braid group $\PB_n$.  Brendle and the author of this article \cite{level4} proved that $\B_n[4]$ is equal to the group generated by squares of Dehn twists in $D_n$ and that both of these are equal to the group $\PB_n^2$.  

A'Campo proved that for $n=2g+1$ the image of $\PB_n$ in $\Sp_{2g}(\Z)$ is the level two congruence subgroup $\Sp_{2g}(\Z)[2]$.  It follows that for $m$ even the image of $\B_n[m]$ in $\Sp_{2g}(\Z)$ is $\Sp_{2g}(\Z)[m]$.

\begin{problem}
Describe the image of $\B_n[m]$ in $\Sp_{2g}(\Z)$ for $m$ odd.
\end{problem}

The group $\PB_n = \B_n[2]$ is generated by Dehn twists and the group $\B_n[4]$ is generated by squares of Dehn twists \cite{level4}.

\begin{problem}\label{bgen}
Find natural generating sets for the groups $\B_n[m]$.
\end{problem}

One generating set for the case of $m$ an odd prime was given by Styliankis \cite{haris}.  The answer to Problem~\ref{bgen} will be more complicated in general than for the cases of $\B_n[2]$ and $\B_n[4]$.  Those are the only two that contain all squares of Dehn twists.  Since $\BI_n$ is generated by squares of Dehn twists about odd curves (as above), the group $\BI_n$ is automatically contained in the stated generating set.

\begin{problem}
\label{kordek}
Compute the abelianizations (and other homology groups) of the groups $\B_n[m]$.
\end{problem}

Kordek points out the following application: if the ranks of the abelianizations of $\B_n[2m]$ are unbounded for fixed $n$ and varying $m$, then $\BI_n$ is not finitely generated (in fact its abelianization has infinite rank).  The proof of this uses the fact, due to Ka\v zdan \cite{dak}, that the congruence subgroups $\Sp_{2g}(\Z)[2m]$ have trivial first integral cohomology. 

Styliankis \cite{haris} has shown that for $p$ an odd prime the quotient group $\B_n[p]/\B_n[2p]$ is isomorphic to the symmetric group $\Sigma_n$.  This in particular implies that the cohomology of $\B_n[2p]$ is a representation of the symmetric group.

In forthcoming work, the first rational homology of $\B_n[4]$ is computed by Kordek and the author.


\section{Generating with torsion}

Almost 50 years ago, MacLachlan proved that $\Mod(S_g)$ is generated by torsion \cite{modulus}.  In the time since, there have been many (successful) attempts to sharpen this result: Which orders are possible?  What is the smallest number of generators needed?

Most of the work in this direction has to do with involutions, elements of order two.  McCarthy and Papadopoulos proved that $\Mod(S_g)$ is generated by involutions, and in fact it is normally generated by a single involution \cite{mcpapinv}.  Luo showed that $\Mod(S_g)$ is generated by $12g+2$ involutions \cite{luo}.  Brendle and Farb sharpened this further: for $g \geq 3$ only 6 involutions are needed, independent of $g$.  In Kassabov \cite{kassabov} goes even further: only four involutions are needed if $g \geq 7$.

We know that $\Mod(S_g)$ cannot be generated by two involutions, because any such group is virtually cyclic.  So we are left with a tantalizing problem, suggested by Brendle--Farb and Kassabov.

\begin{question}
\label{k}
Is $\Mod(S_g)$ generated by three involutions or not?
\end{question}

Of course there are other orders besides two.  Monden \cite{monden} showed for $g \geq 3$ that $\Mod(S_g)$ is generated by three elements of order 3 and also by four elements of order 4.  Yoshihara \cite{yoshihara} has recently shown that $\Mod(S_g)$ is generated by four elements of order 6 when $g \geq 5$.

What about order 5?  One obstacle here is that $\Mod(S_3)$ does not even have an element of order 5.  In fact there are subtle number theoretic conditions that determine whether $\Mod(S_g)$ contains an element of order $k$.  So it seems at first that there is no hope of extending Monden's results to periodic elements of higher periods.  

However, it is a theorem of Tucker \cite{tucker} that for any $k \geq 2$, the group $\Mod(S_g)$ has elements of order $k$ as long as $g \geq (k-1)(k-3)$ (in fact his result is stronger). Building on this, Lanier \cite{lanier} proved that for any $k \geq 6$ the group $\Mod(S_g)$ is generated by three elements of order $k$ when $g \geq (k-1)^2+1$; the new phenomenon here is that the number of generators is independent of both $g$ as well as $k$.  Moreover, Lanier's elements are all conjugate to each other.

Given Tucker's result, we can ask the analogue of Question~\ref{k}.

\begin{question}
\label{x}
For which $k$ can $\Mod(S_g)$ be generated by two elements of order $k$?
\end{question}

The maximum order of a periodic element of $\Mod(S_g)$ is $4g+2$, and Korkmaz showed that $\Mod(S_g)$ is generated by only two of these \cite{kork}.  Here the order $4g+2$ varies with $g$, so this result is in a slightly different direction than Question~\ref{x}

Lanier chooses his periodic elements very carefully: they are all realized by rotations of $S_g$ in $\R^3$.  What about other periodic elements?  With Lanier, the author  proved that every nontrivial periodic element of $\Mod(S_g)$---except for the hyperelliptic involution---is a normal generator for $\Mod(S_g)$.  The following question is asked in our paper \cite{lanierm}.

\begin{problem}
Given a periodic element of $\Mod(S_g)$, how many of its conjugates are needed to generate $\Mod(S_g)$?  Can this number be taken to be independent of the element chosen?
\end{problem}


\section{Generators and relations for Torelli groups}

As mentioned, the Torelli group $\I(S_g)$ is the kernel of the representation $\Mod(S_g) \to \Sp_{2g}(\Z)$ given by the action of $\Mod(S_g)$ on $H_1(S_g;\Z)$.  We may think of the Torelli group as being the non-linear, more mysterious part of the mapping class group.

The Torelli group arises in algebraic geometry because it is the fundamental group of Torelli space, the space of Riemann surfaces with homology framings.  This space is an infinite regular cover of moduli space and it is the natural domain of the period mapping, which maps a Riemann surface with homology framing to Siegel space; see, e.g. Mess' paper \cite{Mess}.

Topologists became interested in the Torelli group because of the following fact: every homology 3-sphere can be obtained from $S^3$ by cutting along a surface and regluing by the product of the original gluing map and an element of the Torelli group.  In his survey paper, Johnson credits the resurgence of interest in the Torelli group to Birman \cite{JohnsonSurvey}.

In 1983 Johnson \cite{Johnson1} proved that $\I(S_g)$ is finitely generated by bounding pair maps for $g \geq 3$.  Mess proved that $\I(S_2)$ is isomorphic to a free group of infinite rank \cite{Mess}.  The following is the most fundamental open question about the Torelli group.  It is asked by Birman in her book \cite[Appendix, Problem 29]{jsb}.

\begin{question}
\label{torellifp}
Let $g \geq 3$.  Is the Torelli group $\I(S_g)$ finitely presented?
\end{question}

The Torelli group has been studied for a century, going back to the work of Dehn, Nielsen, and Magnus.  Still this question is wide open.  Putman has found an infinite presentation for $\I(S_g)$ where all the relations come from geometric considerations \cite{PutmanPresentation}, but there is little hope of turning this into a finite presentation.  Morita and Penner \cite{mp} give a sort of presentation for $\I(S_g)$.  Their presentation has a geometric origin, but the generating set consists of infinitely many copies of each element of the group and the list of relations contains all relations between all elements (and indeed infinitely many copies of each relation).

\subsection*{Generating sets} There is still more to be said about the finite generation of $\I(S_g)$.  The number of elements in Johnson's original finite generating set for $\I(S_g)$ is exponential in $g$.  Johnson also proved \cite{Johnson3} that the abelianization of $\I(S_g)$ has rank $2g+1 \choose 3$, and he conjectured \cite{Johnson1,JohnsonSurvey} that there should be a generating set for $\I(S_g)$ whose size is cubic in $g$.  This conjecture was proved by Putman, who found a generating set for $\I(S_g)$ with $57 {g \choose 3}+2g+1$ elements \cite{putmansmall}.  While this is cubic in $g$, it leaves open the following question, asked explicitly by Johnson \cite{Johnson1}.

\begin{question}
\label{min}
Is there a generating set for $\I(S_g)$ with exactly $2g +1 \choose 3$ elements?
\end{question}
Johnson whittled down his original generating set for $\I(S_3)$ to one with ${7 \choose 3}=35$ elements, thus answering his own question affirmatively in this case.

Johnson's proof that $\I(S_g)$ requires some difficult computations.  His arguments essentially reduce the problem to the case $g=3$, but all of the difficult computations are still required even for this case.  Putman's paper \cite{putmansmall} also gives a conceptual reduction to the case $g=3$.  Still we have the following problem.

\begin{problem}
Find a (more) conceptual proof that $\I(S_3)$ is finitely generated.
\end{problem}

One approach to the previous problem would be to simply elucidate Johnson's proof, for instance by making his computations more transparent (say, realizing them as relations among push maps).  A totally different approach would be to use the fact that the normal closure of any hyperelliptic involution contains $\I(S_g)$ as a subgroup of index 2; see, e.g., the author's paper with Lanier \cite{lanierm}.

Even though $\I(S_2)$ is infinitely generated, there is still a basic open question about its generation.  Mess shows that $\I(S_2)$ is generated by Dehn twists, one for each symplectic splitting of $H_1(S_2;\Z)$.  He posed the following \cite{Mess}.

\begin{problem}
\label{g2}
Give an explicit infinite generating set for $\I(S_2)$.
\end{problem}
Mess conjectured in his paper that one can construct a generating set by choosing a hyperbolic metric on $S_2$ and choosing for each symplectic splitting the Dehn twist about the shortest curve inducing that splitting.  

We expect that one can solve Problem~\ref{g2} using the action of $\I(S_2)$ on the complex of cycles $\mathcal{B}(S_2)$ of Bestvina, Bux, and the author \cite{BBM}.
 
\subsection*{Two-generator subgroups} Related to presentations of the Torelli group we also have the following question.  Our inability to answer it underscores our lack of knowledge about the basic algebraic structure of the Torelli group.

\begin{question}
\label{2elts}
Let $g \geq 3$.  Suppose that $f$ and $h$ are elements of $\I(S_g)$.  Is it true that $f$ and $h$ either commute or generate a free group?
\end{question}

With Leininger, the author answered the analogous question for pure braid groups in the affirmative \cite{LM}.  The pure braid group is the subgroup of the braid group acting trivially on the homology of the punctured disk, and so we may think of the pure braid group as one analogue of the Torelli group (probably a better analogue is the braid Torelli group; see \cite{BMP}).  It is also interesting to consider the analog of Question~\ref{2elts} for pure surface braid groups or for the Johnson kernel.

Question~\ref{2elts} is true when $f$ and $h$ are Dehn twists (cf. \cite[Theorem 3.14]{primer}).  However it is open when $f$ and $h$ are both bounding pair maps, or when one is a bounding pair map and one is a Dehn twist.

Leininger has pointed out that when $f$ and $h$ are bounding pair maps with the property that each of the four curves intersects each of the two curves of the other bounding pair in a single point, the group generated by $f$ and $h$ is free.  This is because the two bounding pairs lie inside a torus with two boundary components, and the corresponding bounding pair maps become trivial after one of the boundary components is filled in.  In other words, the two bounding pair maps lie in a subgroup isomorphic to $\pi_1(S_{1,1}) \cong F_2$.

\subsection*{Homology of Torelli groups} Beyond Question~\ref{torellifp} we also have the following basic question, discussed, for example, by Farb \cite[Problem 5.11]{FarbProblems}.

\begin{question}
\label{torellihom}
For which natural numbers $g$ and $k$ is $H_k(\I(S_g);\Z)$ finitely generated?
\end{question}

By Mess' theorem, we have a complete answer for $g \leq 2$, and by Johnson's theorem we have a complete answer for $k=1$.  With Bestvina and Bux, the author showed \cite{BBM} that the cohomological dimension of $\I(S_g)$ is $3g-5$ (the lower bound was previously given by Mess \cite{MessUnit}).  We also showed that $H_{3g-5}(\I(S_g);\Z)$ is not finitely generated, sharpening an earlier result of Akita \cite{Akita}.  It is a theorem of Johnson--Millson--Mess that $H_3(\I_3;\Z)$ is not finitely generated \cite{Mess}.

Very recently, Gaifullin \cite{aag} has proved that the groups $H_k(\I(S_g);\Z)$ are infinitely generated for $2g-3 \leq k \leq 3g-6$.  It was explicitly asked in our paper with Bestvina and Bux if this was true \cite[Question 8.11]{BBM}.

If a group is finitely presented then its second integral homology group is finitely generated, and so the case of $k=2$ in Question~\ref{torellihom} is of particular interest.  

There is a folk conjecture that for fixed $k$ the group $H_k(\I(S_g);\Z)$ is finitely generated once $g$ is large enough; see, e.g. \cite[Conjecture 5.12]{FarbProblems}.  

\subsection*{The Johnson kernel} There are many variants of the Torelli group, all interesting for different reasons.  The most well-studied subgroup of the Torelli group is the Johnson kernel $\K(S_g)$, so called because it is the kernel of the Johnson homomorphism $\tau : \I(S_g) \to \wedge^3 H_1(S_g;\Z)$.  Johnson proved that this kernel is equal to the subgroup of $\I(S_g)$ generated by Dehn twists about separating curves \cite{Johnson2}.  

Most of our questions about $\I(S_g)$ are open for $\K(S_g)$ as well.  It was only very recently shown that $\K(S_g)$ is finitely generated.  This was originally shown by Ershov and He \cite{ershovhe} for $g \geq 12$ and later extended to $g \geq 4$ by Church, Ershov, and Putman \cite{CEP}.  Dim\c{c}a and Papadima had shown a few years earlier that the torsion-free part of the abelianization was finitely generated \cite{DP}.

\begin{problem}
Find an explicit finite generating set for $\K(S_g)$.
\end{problem}

We can also ask about the size of a minimal generating set for $\K(S_g)$, and the rank of $H_1(\K(S_g);\Z)$.  We also have the following question.

\begin{question}
Is $\K(S_3)$ finitely generated or not?
\end{question}

The following natural question is made more intriguing by the Church--Ershov--Putman result.

\begin{question}
Is $\K(S_g)$ finitely presented or not?
\end{question}

Also, it is known \cite{BBM} that the cohomological dimension of $\K(S_g)$ is $2g-3$ but the following is still open for $g \geq 3$.

\begin{question}
Is $H_{2g-3}(\K(S_g);\Z)$ infinitely generated?
\end{question}

We expect that one can use the techniques used for the Torelli group to answer this question; see \cite{BBM}.

\subsection*{Other terms of the Johnson filtration} The Johnson filtration of $\I(S_g)$ is the sequence of groups $\N_k(S_g)$ defined by
\[
\N_k(S_g) = \ker \left(  \Mod(S_g) \to \Out(\pi/\pi^k) \right)
\]
where $\pi = \pi_1(S_g)$ and $\pi^k$ is the $k$th term of its lower central series.  The groups $\N_2(S_g)$ and $\N_3(S_g)$ are equal to $\I(S_g)$ and $\K(S_g)$, respectively.  Magnus proved that the intersection of the $\N_k(S_g)$ (as $k$ varies) is trivial \cite{magnus}.

Church, Ershov, and Putman proved that each term of the Johnson filtration is finitely generated when $k \geq 3$ and $g \geq 2k-1$.

\begin{problem}
Find explicit finite generating sets for the terms of the Johnson filtration.
\end{problem}

It is even an open problem to find simple, explicit, infinite generating sets.  Here is one proposal.  For any $f \in \N_k(S_g)$ and any curve $c$ in $S_g$ the twist difference $[T_c,f] = T_cT_{f(c)}^{-1}$ lies in $\N_k(S_g)$.

\begin{question}
Are the terms of the Johnson filtration generated by twist differences?
\end{question}

This question is already interesting for the Johnson kernel, where it is not known.  We can refine the question by asking if the terms of the Johnson filtration are generated by twist differences on subsurfaces of uniformly small Euler characteristic (that is, the Euler characteristic only depends on $k$).  Church and Putman that there does exist some (infinite) generating set with this property \cite{CP}.

It is known that the cohomological dimension of $\N_k(S_g)$ lies between $g-1$ and $2g-3$ for $k \geq 4$; see \cite{BBM}.  The lower bound comes from an explicit construction of a free abelian subgroup of rank $g-1$.  The upper bound comes from the fact that the cohomological dimension of $\K(S_g)$ is $2g-3$.   Since cohomological dimension is non-increasing under passage to subgroups, it follows that as we increase $k$ (and fix $g$) the cohomological dimension must eventually stabilize.

\begin{problem}
Compute the cohomological dimension of $\N_k(S_g)$ for $k \geq 4$ and $g \geq 3$.
\end{problem}

As a first start we have the following problem.

\begin{problem}
Compute the maximal rank of an abelian subgroup of $\N_k(S_g)$ for $k \geq 4$.
\end{problem}

We conjecture that the maximal rank in the above problem is $g-1$, the lower bound given in Farb's problem list  \cite[Theorem 5.10]{FarbProblems}.  One feature to note is that for $k \geq 4$ the groups $\N_k(S_g)$ contain no multitwists \cite[Theorem A.1]{BBM}.

\subsection*{The Magnus filtration} There is a natural variant of the Johnson filtration, called the Magnus filtration.  These groups are obtained by replacing $\pi^k$ in the definition of the Johnson filtration with the $k$th term of the derived series of $\pi_1(S_g)$.  This filtration was defined by McNeill \cite{taylor} who showed, among other things, that the the quotient of each term by the next is infinitely generated.  Beyond her paper, these groups have been unexplored.

\begin{problem}
Investigate the basic properties of the terms of the Magnus filtration.  For instance what is a simple, explicit generating set?  What are the abelianizations?
\end{problem}


\section{Virtual surjection onto the integers}

The mapping class group $\Mod(S_g)$ is perfect for $g \geq 3$, meaning that it is equal to its own commutator subgroup.  For $g \leq 2$ the commutator subgroup of $\Mod(S_g)$ has finite index.  So for every $g$ all homomorphisms $\Mod(S_g) \to \Z$ are trivial.  The following is one of the most fundamental open questions about the mapping class group.

\begin{question}
\label{vfb}
Is there a finite-index subgroup of $\Mod(S_g)$ that surjects onto $\Z$?
\end{question}

Since $\Mod(S_1)$ is virtually free the answer is `yes' for $g=1$.  Taherkhani \cite{t} and McCarthy \cite{mc} proved that the answer is also `yes' for $g=2$.  Thus, we may restrict attention to $g \geq 3$.

McCarthy \cite{mc} and Hain \cite{hain} answered Question~\ref{vfb} in the negative for all subgroups of finite index containing the Torelli group, and Putman \cite{andynote} extended Hain's result to all subgroups of finite index containing the Johnson kernel.  Finally, Ershov and He extended these results to all terms of the Johnson filtration \cite{ershovhe}.

\begin{problem}
Compute the abelianizations of the various Morita subgroups, congruence subgroups, and liftable subgroups of $\Mod(S_g)$ defined in Section 2.  
\end{problem}

Mortia's subgroups contain the Johnson kernel, and so by Putman's result these abelianizations are torsion.

Putman and Wieland \cite{pw} conjecture that the answer to Question~\ref{vfb} is `no.'  They further prove that if their conjecture holds for $g=3$ then it holds for all $g \geq 3$. 

Putman and Wieland also relate Question~\ref{vfb} to a question about the so-called Prym representations of $\Mod(S_g)$.  In short, a Prym representation is the linear representation given by the action of $\Mod(S_g)$ on the first homology of some finite, characteristic cover (in order for this representation to be well-defined, we must fix a basepoint in $S_g$; see their paper for details).  Prym representations have been studied by Looijenga \cite{looi}, Grunewald--Lubotzky \cite{gl}, and Grunewald--Larsen--Lubotzky--Malestein \cite{gllm}.  Putman and Wieland conjecture the following.

\begin{conjecture}
\label{pw}
For $g \geq 2$ all nonzero orbits of all Prym representations of $\Mod(S_g)$ are infinite.
\end{conjecture}

Putman and Wieland proved that this conjecture is essentially equivalent to their conjecture that the answer to Question~\ref{vfb} is `no.'  In particular, to answer Question~\ref{vfb} in the negative, it is enough to prove Conjecture~\ref{pw} for the case of a surface of genus two with a single boundary component.

Farb and Hensel proved that the analog of Conjecture~\ref{pw} for graphs is true \cite{fh}.


\section{Ivanov's Metaconjecture}
\label{sec:normal}

The problems in the next two sections are motivated by the following general question:
\begin{quote}
\emph{What can we say about the typical (or general) normal subgroup of the mapping class group?}
\end{quote}

In this section we discuss the automorphism groups of normal subgroups.  The story begins with a theorem of Ivanov \cite{ivanov}, which states that
\[
\Aut \Mod(S_g) \cong \Mod^\pm(S_g)
\]
for $g \geq 3$.  Here $\Mod^\pm(S_g)$ is the extended mapping class group, the group of homotopy classes of all homeomorphisms of $S_g$ (including the orientation-reversing ones).  

Ivanov's theorem gives the isomorphism explicitly: the natural map $\Mod^\pm(S_g) \to \Aut \Mod(S_g)$ given by conjugation is an isomorphism.  A key ingredient is the accompanying theorem that the automorphism group of the complex of curves $\C(S_g)$ is also isomorphic to $\Mod^\pm(S_g)$.  To compute $\Aut \Mod(S_g)$, Ivanov first shows that an automorphism must preserve powers of Dehn twists, and then uses this to reduce to his theorem about $\C(S_g)$.

After Ivanov's work, many other similar theorems were proved, each saying that the automorphism group of some normal subgroup of $\Mod(S_g)$ or some simplicial complex associated to $S_g$ is isomorphic to $\Mod^\pm(S_g)$.  In response to many of these results, Ivanov posed the following \cite{15}.

\begin{meta}
Every object naturally associated to a surface $S$ and having a sufficiently rich structure has $\Mod^\pm(S_g)$ as its group of automorphisms.  Moreover, this can be proved by a reduction to the theorem about the automorphisms of $\C(S)$.
\end{meta}

One type of ``objects'' to consider are the normal subgroups of $\Mod(S_g)$.  Farb and Ivanov \cite{farbivanov} proved that the Torelli group $\I(S_g)$ satisfies the metaconjecture.  Brendle and the author did the same for the Johnson kernel \cite{kg,kgadd,kgerr}.  Bridson, Pettet, and Souto \cite{bps} then announced a vast generalization: all terms of the Johnson filtration have automorphism group isomorphic to $\Mod^\pm(S_g)$ for $g \geq 4$.  It is rather startling that these arbitrarily small subgroups of the mapping class group somehow remember the structure of the original group.  (In fact much more is true: in each case it was further shown that the abstract commensurator group of the given subgroup is $\Mod^\pm(S_g)$.)  

In the spirit of Ivanov, each of the above results is proved by showing that some appropriate combinatorial complex has automorphism group isomorphic to $\Mod^\pm(S_g)$: the complex of curves, the complex of separating curves and bounding pairs, the complex of separating curves, and the complex of shirts and straitjackets.  In fact, there were many other papers proving that the automorphism group of some combinatorial complex associated to $S_g$ is isomorphic to $\Mod^\pm(S_g)$.  See the author's survey with Brendle \cite{survey} or the paper by McCarthy--Papadopoulos \cite{mcpap} for comprehensive surveys.

The result of Bridson--Pettet--Souto gives many examples of normal subgroups satisfying Ivanov's metaconjecture.  One might hope that all nontrivial normal subgroups of $\Mod(S_g)$ satisfy the metaconjecture.  However, this is certainly not the case.  Dahmani--Guirardel--Osin \cite{dgo} constructed examples of infinitely generated, free, normal subgroups of $\Mod(S_g)$; the automorphism groups of these groups are uncountable, and hence definitely not isomorphic to $\Mod^\pm(S_g)$.

\subsection*{The metaconjecture for normal subgroups} A recent theorem of the author with Brendle gives a large class of normal subgroups of $\Mod(S_g)$ that satisfy Ivanov's metaconjecture \cite{meta}.  We will state a simplified  version here.

Say that a subsurface of $S_g$ is \emph{small} if it is contained as a non-peripheral subsurface of a subsurface of $S_g$ of genus $k$ with connected boundary where $k = \lfloor g/3 \rfloor$.  An element of $\Mod(S_g)$ has \emph{small support} if it has a representative whose support is small.

The theorem we proved is that if $N$ is a normal subgroup of $\Mod(S_g)$ that contains a nontrivial element of small support, then
\[
\Aut N \cong \Mod^\pm(S_g).
\]
Additionally, we prove that the abstract commensurator group of $N$ is $\Mod^\pm(S_g)$.  The examples of Dahmani--Guirardel--Osin are all pseudo-Anosov, which means that the support of every nontrivial element is the whole surface $S_g$.  Our hypothesis exactly rules out these types of examples.

Our theorem implies the result announced by Bridson--Pettet--Souto about the terms of the Johnson filtration (only for $g \geq 7$).  It also applies to many other normal subgroups, including the terms of the Magnus filtration, discussed above.  Using our theorem one can define many other normal subgroups satisfying Ivanov's metaconjecture, for instance, the group generated by 100th powers of Dehn twists about separating curves of genus 7 (our result applies for $g \geq 25$).

In the statement of the theorem there is a $g/3$ that should morally be $g/2$; that is, half the surface instead of one-third.  Say that a nonseparating subsurface of $S_g$ is \emph{not large} if it is homeomorphic to a proper subsurface of its complement.  And say that an arbitrary subsurface is not large if it is contained in a nonseparating subsurface that is not large.  With Brendle, we make the following conjecture \cite{meta}.

\begin{conjecture}
\label{large}
Suppose that $N$ is a normal subgroup of $\Mod(S_g)$ that contains an element whose support is not large.  Then the automorphism group (and also the abstract commensurator group) of $N$ is isomorphic to $\Mod^\pm(S_g)$.
\end{conjecture}

Chen has pointed out that whenever we have a normal subgroup $N$ of $\Mod(S_g)$ with $\Aut N \cong \Mod^\pm(S_g)$ then every normal subgroup of $\Mod(S_g)$ isomorphic to $N$ is equal to $N$; see \cite[Corollary 1.2]{meta}.  We can ask if the same holds when $N$ is not assumed to be normal.

\begin{question}
Is it true that if $N$ is a normal subgroup of $\Mod(S_g)$ with $\Aut N \cong \Mod^\pm(S_g)$ then every subgroup of $\Mod(S_g)$ isomorphic to $N$ is equal to $N$?  Is this true for the Torelli group $\I(S_g)$?
\end{question}

\subsection*{The metaconjecture for curve complexes} Our theorem with Brendle is proved by showing that there is a large class of simplicial complexes associated to $S_g$ whose members all have automorphism group isomorphic to $\Mod^\pm(S_g)$.  Our complexes have vertices corresponding to proper, connected subsurfaces of $S_g$ and edges for disjointness.  The hypotheses of our theorem require that the simplicial complexes are connected, that they have a small element, and that they admit no exchange automorphisms, which are automorphisms that interchange two vertices and fix all other vertices.  There are a number of ways in which one might try to improve on our result.

\begin{problem}
\label{metagen}
Improve on our theorem with Brendle about automorphisms of simplicial complexes associated to $S_g$ in the following ways:
\begin{enumerate}
\item replace ``small'' with ``not large,''
\item allow vertices corresponding to disconnected subsurfaces, and
\item allow edges for configurations other than disjointness.
\end{enumerate}
\end{problem}
The last two items would allow for complexes such as the pants complex, whose vertices correspond to disconnected subsurfaces and whose edges do not correspond to disjointness.  

In our paper \cite{meta} we conjecture that there is a stronger statement than the one suggested by item (1) in Problem~\ref{metagen}.  Specifically, we conjecture that the small assumption can be removed altogether, if we keep the assumption of connectivity.

\begin{conjecture}
Let $\C_A(S_g)$ be a complex of regions that is connected and admits no exchange automorphisms.  Then the natural map
\[ \Mod^\pm(S_g) \to \Aut \C_A(S_g) \]
is an isomorphism.
\end{conjecture}

\subsection*{Beyond the metaconjecture} With Brendle, we also give the following generalization of Ivanov's metaconjecture \cite{survey}.

\begin{genmeta}
Suppose that $Y$ is a nice space.  Every object naturally associated to $Y$ and having a sufficiently rich structure has $\Out \pi_1(Y)$ as its group of automorphisms. 
\end{genmeta}

This is indeed a generalization of Ivanov's metaconjecture because of the Dehn--Nielsen--Baer theorem \cite[Theorem 8.1]{primer} which states that $\Out \pi_1(S_g) \cong \Mod^\pm(S_g)$.  One example of a nice space $Y$ is the connected sum of $n$ copies of $S^2 \times S^1$.  The fundamental group of this manifold is the free group $F_n$. Scott \cite{scott} has shown that many complexes associated to this manifold have automorphism group $\Out(F_n)$.  Other specific complexes associated to this manifold have been studied by Aramayona--Souto \cite{as}, Bestvina--Bridson \cite{bestb}, and Pandit \cite{sp,sp2}.


\section{Normal right-angled Artin subgroups}

As mentioned, Dahmani--Guirardel--Osin constructed free normal subgroups in $\Mod(S_g)$.  In their subgroups, all of the nontrivial elements are pseudo-Anosov \cite{dgo}.  Their subgroups are constructed by taking the normal closure of some finite collection of pseudo-Anosov mapping classes.

In joint work with Clay and Mangahas \cite{cmm}, we construct similar free groups with partial pseudo-Anosov elements instead of pseudo-Anosov elements.  By our theorem with Brendle about automorphisms of normal subgroups above, the supports of these partial pseudo-Anosov elements must be large in the sense of the previous section.  In other words the supports take up more than half the surface.

By varying our construction we also produce normal subgroups of $\Mod(S_g)$ of the following form:
\[ \displaystyle\BigFreeProd_\infty \left(F_\infty \times F_\infty\right)\, \quad  \text{ and } \quad \displaystyle\BigFreeProd_\infty \left(F_\infty \times F_\infty \times \Z\right).\]
For surfaces with punctures, we also construct normal subgroups isomorphic to \[ \displaystyle\BigFreeProd_\infty \left(F_\infty \times \Z\right). \]

To create a normal subgroup of the first type in $\Mod(S_g)$, we take $g$ to be even and take the normal closure of a (high power of a) partial pseudo-Anosov element $f$ of $\Mod(S_g)$ whose support is homeomorphic to $S_{g/2}^1$.  The element $f$ has a conjugate in $\Mod(S_g)$ supported in the complement; this is how the commuting elements arise.  To construct an example of the second type of group, we take the normal closure of the same $f$ and a power of $T_c$, where $c$ is a curve in $S_g$ of genus $g/2$ (that is, $c$ cuts $S_g$ in half).  

\begin{conjecture}
Any normal subgroup of $\Mod(S_g)$ isomorphic to a right-angled Artin group is isomorphic to a free product of groups from the following list:
\[
F_\infty\, , \quad \displaystyle\BigFreeProd_\infty \left(F_\infty \times F_\infty\right)\, , \quad  \text{ and } \quad \displaystyle\BigFreeProd_\infty \left(F_\infty \times F_\infty \times \Z\right).
\]
\end{conjecture}

We also conjecture that these are the only examples of normal subgroups that do not have automorphism group the extended mapping class group.

\begin{conjecture}
A normal subgroup of $\Mod(S_g)$ either has automorphism group $\Mod^{\pm}(S_g)$ or it is isomorphic to a right-angled Artin group.
\end{conjecture}


\section{Cohomology of the mapping class group}

In this section we discuss the following central open problem. Throughout we take our coefficients for homology to be $\Q$.

\begin{problem}
\label{h}
Compute $H^k(\Mod(S_g))$ for each $g$ and $k$.
\end{problem}

A more modest question is to determine which $H^k(\Mod(S_g))$ are nonzero.

One of the reasons the cohomology of the mapping class group is important because it plays the same role for surface bundles that the classical characteristic classes play for vector bundles.

In this vein, we may think of an element of $H^k(\Mod(S_g))$ as follows.  It is a function $\alpha$ that assigns to each $S_g$-bundle over a space $B$ an element of $H^k(B)$; the function $\alpha$ has to satisfy several naturality properties, the most important being that it behaves naturally under pullbacks: if $E$ is a bundle over $B$ and $f : A \to B$ then $\alpha(f^\star(E)) = f^\star(\alpha(E))$.  By the naturality of pullback, an element of $H^k(\Mod(S_g))$ is completely determined by its behavior on $S_g$-bundles over $k$-manifolds, in which case $H^k(B)$ is a number.

\subsection*{Background} We begin by describing what is known about the cohomology of $\Mod(S_g)$.  Throughout, we refer to Figure~\ref{dark}.

\begin{figure}
\includegraphics[width=\textwidth]{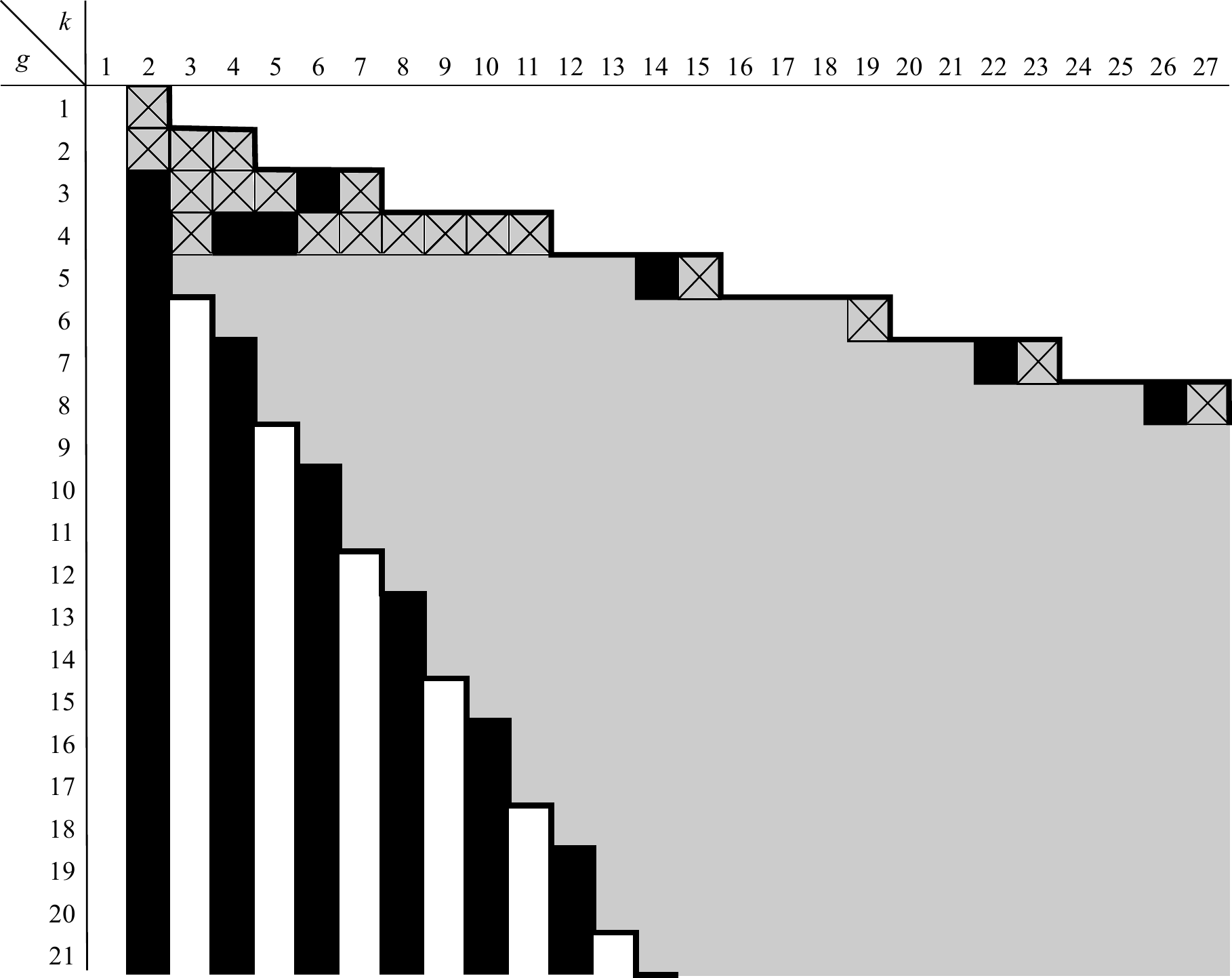}
\label{dark}
\caption{A summary of what is known about the cohomology of $\Mod(S_g)$; filled boxes denote nontrivial cohomology groups, the shaded region denotes the unstable region, and in this region X's denote trivial cohomology groups}
\end{figure}

In the 1980s Harer proved that the virtual cohomological dimension of $\Mod(S_g)$ is $4g-5$ for $g \geq 2$, and so $H^k(\Mod(S_g))$ is trivial for $k > 4g-5$.  The latter fact was also proved by Culler--Vogtmann.  This explains the empty region at the top right of Figure~\ref{dark}: all of these cohomology groups are zero.

Harer proved that the mapping class groups satisfy homological stability \cite{harerstab}.  This means that for fixed $k$, and $g$ large enough, the groups $H^k(\Mod(S_g))$ do not depend on $g$.  Specifically, this has been shown to hold in the ``stable range'' $k \leq 2 \lfloor \frac{g}{3} \rfloor$, by the work of Harer \cite{harerimp}, Ivanov \cite{ivstab1}, Wahl \cite{wahl}, Randal-Williams \cite{rw}, and Boldsen \cite{boldsen}. 

Mumford conjectured that in the stable range the cohomology is generated by the Morita--Mumford--Miller classes $\kappa_1,\kappa_2,\dots$.  Each $\kappa_i$ is an element of $H^{2i}(\Mod(S_g))$.  The $\kappa_i$ are also called  tautological classes.

Here is how $\kappa_1$ is defined.  As above, $\kappa_1$ should assign to an $S_g$-bundle $E$ over a surface $S_h$ an element of $H^2(S_h)$ (that is, a rational number).  The bundle $E$ has a vertical tangent bundle, given by the tangent planes to the fibers.  This 2-plane bundle has an Euler class $e \in H^2(E)$.  The class $\kappa_1$ is the image of $e^2$ under the Gysin map $H^2(E) \to H^2(S_h)$.  In other words, $\kappa_1(E)$ is obtained as the intersection number of three generic sections of the vertical tangent bundle over $E$.  In general, the classes $\kappa_i$ are defined as the images under the Gysin homomorphism of $e^{i+1}$ (so intersect $i+2$ sections).

Mumford's conjecture was proved in 2006 by Madsen--Weiss \cite{mw}; see also \cite{egm,hatcher,wahl2}.  This gives a complete picture of the cohomology in the stable range: it is the polynomial algebra generated by the $\kappa_i$ (for fixed $g$ these are necessarily trivial for $2i > 4g-5$).  In Figure~\ref{dark} this stable range is the region on the bottom left; the reason for the vertical stripes is Harer's stability theorem, plus the fact that the $\kappa_i$ all have even degree.

In order to solve Problem~\ref{dark} we should thus focus our attention on the unstable range, shaded in gray in Figure~\ref{dark}.  A few of the cohomology groups here are known.  In fact, we have complete computations of $H^\ast(\Mod(S_g))$ for $g \leq 4$. The black boxes in the first four rows of Figure~\ref{dark} represent 1-dimensional cohomology groups and the X's represent trivial groups.  The cohomology groups for $g$ equal to 2, 3, and 4 were computed by Igusa in 1960, Looijenga in 1991, and Tommasi in 2005, respectively \cite{igusa,edloo,tommasi}.

\subsection*{Dark matter} While very few unstable classes are known explicitly, we do know that many other classes exist.  Harer and Zagier \cite{hz} computed the Euler characteristic of moduli space $\mathcal{M}_g$, thought of as an orbifold.  They showed that
\[
\chi(\mathcal{M}_g) = \frac{\zeta(1-2g)}{2-2g} \sim (-1)^g \frac{(2g-1)!}{(2-2g)2^{2g-1}\pi^{2g}}.
\]
In particular, this number grows super-exponentially in $g$.  This implies that the Betti numbers of torsion-free subgroups of finite index in $\Mod(S_g)$ grow super-exponentially in $g$.  Harer and Zagier promote this to the statement that the Betti numbers of $\Mod(S_g)$ itself grow super-exponentially in $g$.

The number of stable classes is polynomial in $g$ (since it is a polynomial algebra).  This means that there is a super-exponential number of classes in the unstable range waiting to be discovered.  What is more, since $\chi(\mathcal{M}_g)$ is negative for $g$ even, we know that there are super-exponentially many classes in odd degree.

\begin{problem}
\label{dark}
Find some (or even one!) element of $H^k(\Mod(S_g))$ outside of the stable range.
\end{problem}

To be even more greedy, we would like geometrically meaningful interpretations of the unstable elements.  The first tautological class $\kappa_1$, for example, measures the signatures (and first Pontryagin classes) of surface bundles and relates to the Weil--Petersson volume form on Teichm\"uller space.  Salter \cite{salter} gives similar characterizations of the other $\kappa_i$.  We would like to have geometric interpretations for the unstable classes.  For instance, what does Looijenga's unstable class in $H^6(\Mod(S_3))$ tell us about $S_3$-bundles over 6-manifolds?

Farb has called Problem~\ref{dark} the ``dark matter problem.''  (Of course an important difference from the original dark matter problem is that we know the dark matter exists!)

Church--Farb--Putman \cite{CFP} and Morita--Sakasai--Suzuki independently proved that $H^{4g-5}(\Mod(S_g))$ is trivial, so there are no unstable classes in the top degree.  This might seem like a contradiction at first, since the virtual cohomological dimension of $\Mod(S_g)$ is $4g-5$.  But the virtual cohomological dimension being $4g-5$ only requires that $\Mod(S_g)$ has nontrivial cohomology in degree $4g-5$ with \emph{some} (possibly twisted) coefficients.  This theorem is indicated by several X's in Figure~\ref{dark}.

Very recently, Chan--Galatius--Payne \cite{cgp} proved that $H^{4g-6}(\Mod(S_g))$ is nontrivial for $g = 5$ and $g \geq 7$.  Even more, they show that the dimension grows exponentially in $g$.  This is the first partial answer to Problem~\ref{dark} that fills in infinitely many boxes in Figure~\ref{dark}.

Even if we knew what all of the unstable elements of $H^\ast(\Mod(S_g))$ were, we would still have the problem of understanding the ring structures for each $g$.   In other words, what are the relations? The Faber conjectures \cite{faber} describe the structure of the tautological ring, that is, the ring generated by the tautological classes $\kappa_i$.  The conjectures give explicit proportionality relations between the tautological classes and assert that the tautological ring  for $\Mod(S_g)$ ``looks like'' the algebraic cohomology ring of a smooth projective variety of dimension $g-2$.  Much of the conjectures have been proven; see the surveys by Faber \cite{faber} and Morita \cite{moritaut}.

\subsection*{Torelli groups} Morita proved that the odd $\kappa_i$ vanish on the Torelli group $\I(S_g)$.  He then asked the following \cite[Problem 2]{morita}.

\begin{question}
Are the restrictions of the even $\kappa_i$ to $\I(S_g)$ trivial or not?
\end{question}

Even for $i=2$ the answer is unknown.  In this case, Morita's question asks whether there is an $S_g$-bundle over a 4-manifold, with monodromy in $\I(S_g)$, so that if we intersect four generic sections of the vertical tangent bundle we obtain a nonzero number.


\section{Pseudo-Anosov theory}

The Nielsen--Thurston classification theorem says that mapping classes fall into three types: periodic, reducible, and pseudo-Anosov.  To each pseudo-Anosov element there is an associated stretch factor. 

\subsection*{Genericity} One of the most famous open problems in this direction is the following.

\begin{problem}
Show that pseudo-Anosov mapping classes are generic.  
\end{problem} 

This problem has been answered in one sense.  Maher has shown that a random walk in the mapping class group gives rise to a pseudo-Anosov element with asymptotic probability one \cite{maher}.  Rivin proved the analogous result for the Torelli group \cite{rivin}.

Another notion of genericity is that the proportion of pseudo-Anosov elements in the ball of radius $r$ in a Cayley graph for $\Mod(S_g)$ should be asymptotically one.  This version is open.  Recently  Cumplido and Wiest have shown that this proportion is asymptotically positive \cite{cw}.  

\subsection*{Fried's question}

Another fundamental question is the following.

\begin{question}
\label{real}
Which real numbers are stretch factors of pseudo-Anosov mapping classes?
\end{question}

There are many variants of this question.  The most famous version is due to Fried, who posed the following question \cite[Problem 2]{df}.  In the statement, a bi-Perron unit is an algebraic integer $\alpha$ so that all its algebraic conjugates are either equal to $1/\alpha$ or have absolute value in $(1/\alpha,\alpha)$ (the number $1/\alpha$ is not required to be an algebraic conjugate).

\begin{question}
\label{fried}
Is is true that every bi-Perron unit has a power that is a stretch factor of a pseudo-Anosov mapping class?
\end{question}

Very little is known about this question.  Recently Pankau has shown that every Salem number has a power that is a stretch factor of a pseudo-Anosov mapping class.  A Salem number is a real number $\alpha>1$ so that $1/\alpha$ is an algebraic conjugate and so that all other algebraic conjugates lie on the unit circle.

We remark that the power in Question~\ref{fried} may not be required.

\subsection*{Algebraic degrees} Question~\ref{real} is still (or more?) interesting if we fix the surface we are considering.  Thurston \cite{ThurstonBulletin} proved that the stretch factor of a pseudo-Anosov element of $\Mod(S_g)$ is an algebraic integer with degree between 2 and $6g-6$ (inclusive).  Long showed that if the degree is greater than $3g-3$ then it is even.  (McMullen has given simple proofs of the Thurston and Long restrictions; see the paper by Shin \cite{shin}.)  Strenner then showed that all degrees satisfying the Thurston and Long restrictions do indeed occur \cite{strenner}.  Strenner's proof is not completely constructive.  This leaves open the following.

\begin{problem}
For each $g \geq 2$ and each even $d$ in between $2$ and $6g-6$ and each odd $d$ between $2$ and $3g-3$, construct an explicit pseudo-Anosov element of $\Mod(S_g)$ whose stretch factor has algebraic degree $d$.
\end{problem}

Before Strenner's work, Arnoux--Yoccoz \cite{ay} gave examples of pseudo-Anosov elements of degree $g$.  Also, Shin \cite{shin} gave examples of pseudo-Anosov elements degree $2g$.  

The Arnoux--Yoccoz examples are given in terms of interval exchange transformations; it would be interesting to write them as products of Dehn twists (or other simple mapping classes).\footnote{This problem was recently solved by Liechti--Strenner \cite{ls}.}

Shin found his examples by doing an extensive computer search.  Some algebraic degrees seem to be abundant, and some---such as degree 5 in $\Mod(S_2)$---appear to be rare.  We are thus led to the following question.

\begin{question}
Which algebraic degrees are most common?  Given $g$, is there an algebraic degree that is generic?
\end{question}

We also have the following problem inspired by Strenner's result.

\begin{problem}
Which algebraic degrees appear in the Torelli group $\I(S_g)$?  What about an arbitrary proper normal subgroup of $\Mod(S_g)$?
\end{problem}

\subsection*{Minimal stretch factors} For fixed $g$ the stretch factors in $\Mod(S_g)$ form a discrete subset of the real numbers.  In particular, there is for each $g$ a smallest stretch factor $\lambda_g$.  We thus have the following basic question.

\begin{problem}
\label{smalldil}
Determine the $\lambda_g$.  
\end{problem}

Problem~\ref{smalldil} has only been solved for $g$ equal to 1 and 2.  It is an easy exercise to show that $\lambda_1$ is realized by the conjugacy class of the element of $\Mod(S_1)$ corresponding to the matrix
\[
\begin{bmatrix}  2 & 1 \\ 1 & 1 \end{bmatrix}
\]
And $\lambda_2$ is the largest real root of the polynomial
\[
x^4 -x^3-x^2-x+1
\]
This was proved by Cho--Ham \cite{choham}, who reduced the problem to a finite (but extremely large) check.  Lanneau and Thiffeault \cite{LT} further determined all conjugacy classes in $\Mod(S_2)$ realizing the minimum---there are two---and gave an explicit Dehn twist factorization for each.  Minimal stretch factors have also been found for some braid groups \cite{LT} and some mapping class groups of non-orientable surfaces \cite{bs}.

\begin{question}
What is the smallest stretch factor in $\I(S_g)$?  What about $\K(S_g)$?  Is the former strictly smaller than the latter?
\end{question}

With Farb and Leininger, the author proved \cite{flm2} that, for $g$ fixed, the smallest stretch factor in the $k$th term of the Johnson filtration goes to infinity with $k$.  So we may ask more generally if this sequence is monotonically increasing with $k$. 

\subsection*{Normal closures} With Lanier \cite{lanierm}, the author proved that every pseudo-Anosov mapping class with stretch factor less than $\sqrt{2}$ is a normal generator for $\Mod(S_g)$, that is, its normal closure is the whole group.  Said another way, if a pseudo-Anosov mapping class lies in a proper normal subgroup of $\Mod(S_g)$ then its stretch factor is greater than $\sqrt{2}$.  On the other hand, the only explicit examples of pseudo-Anosov mapping classes that lie in a proper normal subgroup have stretch factor greater than 62; see \cite{flm2}.  So the smallest stretch factor of a pseudo-Anosov mapping class that lies in a proper normal subgroup is between $\sqrt{2}$ and $62$.  

\begin{problem}
Improve the bounds on the smallest stretch factor of a pseudo-Anosov mapping class lying in a proper normal subgroup of the mapping class group.
\end{problem}

With Lanier we also give examples of pseudo-Anosov mapping classes whose normal closure in $\Mod(S_g)$ is equal to any given level $m$ congruence subgroup.  We ask the following.

\begin{question}
Can other normal subgroups, such as the Johnson kernel and the other terms of the Johnson filtration, can be obtained as the normal closure in $\Mod(S_g)$ of a single pseudo-Anosov mapping class?
\end{question}

We also give examples of pseudo-Anosov mapping classes with the property that arbitrarily large powers are normal generators for the mapping class group.  In particular we can take these powers to have arbitrarily large translation lengths on $\C(S_g)$.  The measured foliations associated to our examples have nontrivial symmetry groups.  Since the generic pseudo-Anosov mapping class (in the sense of random walks) has foliations without symmetries \cite{masai}, we arrive at the following problem.

\begin{problem}
Is the generic element of $\Mod(S_g)$ a normal generator, or not?
\end{problem}

The next question is related.  In the statement, the extended Torelli group is the preimage of $\{\pm I\}$ under the standard symplectic representation of $\Mod(S_g)$.  This group is the normal closure of any hyperelliptic involution.

\begin{question}
If the asymptotic translation length of a pseudo-Anosov mapping class is large, can its normal closure be anything other than the mapping class group or a free group?  Can the normal closure be the extended Torelli group?
\end{question}

\subsection*{Small stretch factors} Penner \cite{penner} proved that there are constants $c$ and $C$ so that
\[
\frac{c}{g} \leq \log \lambda_g \leq \frac{C}{g}.
\]
In particular, $\lambda_g \to 1$ as $g \to \infty$.  In particular, while the set of stretch factors in $\Mod(S_g)$ is discrete for fixed $g$, the set of all stretch factors is dense in $(1,\infty)$ (powers of stretch factors are also stretch factors).

We can say that a stretch factor $\lambda$ of an element of $\Mod(S_g)$ is \emph{small} if $g \log \lambda$ is bounded above by Penner's $C$ (or any fixed number greater than or equal to the minimal possible $C$ as above).  Let $P_C$ denote the set of pseudo-Anosov mapping classes $f$ with small stretch factor.  We emphasize here that the set $P_C$ contains elements of $\Mod(S_g)$ for every $g$.  

McMullen gave an alternate proof of Penner's theorem from the point of view of fibered 3-manifolds \cite{ctm}.  Building on work of Fried \cite{fried} and Thurston \cite{norm}, McMullen \cite{ctm} showed that in a single fibered 3-manifold, there are fibers $\Sigma_g$ so that the associated monodromies $f_g : \Sigma_g \to \Sigma_g$ have small stretch factor for $g$ sufficiently large.

Inspired by McMullen's work, we consider the set of all mapping tori $M_f$ where $f \in P_C$. For any such $f$, let $M_f^\circ$ denote the 3-manifold obtained from $M_f$ by removing the suspensions of the singular points of (a representative of) $f$.  With Farb and Leininger \cite{flm}, we proved that the set
\[
\Omega_C = \{ M_f^\circ \mid  f \in P_C \}
\]
is finite.  The removal of the suspensions of singular points is indeed necessary.  Penner \cite{penner} gave examples of elements of $P_C$ with prongs of arbitrarily large degree; but in a given closed 3-manifold, there is an upper bound to the number of prongs.

\begin{question}
What is the smallest nonzero cardinality of $\Omega_C$?  If we vary $C$ in $(1,\infty)$, how does the cardinality grow?
\end{question}

The following related conjecture is due to Farb, Leininger, and the author of this article.

\begin{conjecture}[Symmetry conjecture]
If $f \in \Mod(S_g)$ is a pseudo-Anosov mapping class with stretch factor $\lambda$ satisfying $g \log \lambda < C$, then $f$ is equal to the product of a finite order element and a reducible element, and moreover the absolute value of the Euler characteristic of the support of the reducible element is bounded above by $C$, independently of $g$.
\end{conjecture}

Leininger has shown that if the stretch factor of $f$ is less than $3/2$ then it is the product of a periodic element and a reducible element.  The idea is that in this case we can find a curve $c$ so that $i(c,f(c))$ is at most two, and thus we can find a periodic element taking $f(c)$ back to $c$.  The reducible element constructed here has large complexity.  On the other hand, this is nontrivial progress because there are pseudo-Anosov elements of the mapping class group that cannot be written as the product of a periodic element and a reducible element (apply the work of Bestvina--Fujiwara \cite{bf}).

\subsection*{Pseudo-Anosov surface subgroups} We say that a subgroup of $\Mod(S_g)$ is all pseudo-Anosov if every nontrivial element is pseudo-Anosov.  A version of the following question was raised by Farb and Mosher \cite[Question 1.9]{farbmosher}.

\begin{question}
Does $\Mod(S_g)$ contain an all pseudo-Anosov surface subgroup?
\end{question}

The motivation for the last question is the study of surface bundles over surfaces.  There are no known examples of surface bundles over surfaces that are hyperbolic 4-manifolds.  Of course hyperbolic implies negative curvature in the sense of Gromov.  The latter is known to imply that the monodromy group is all pseudo-Anosov (if there was a surface bundle over a surface with all pseudo-Anosov monodromy that was not hyperbolic, then this would contradict a conjecture of Gromov, which states that any group $G$ with a finite $K(G,1)$ and no Baumslag--Solitar subgroups is hyperbolic).  Therefore, a positive answer to the last question is necessary if there are to exist hyperbolic surface bundles over surfaces.

Leininger--Reid \cite{lr} produced surface subgroups of $\Mod(S_g)$ where every element outside of a single cyclic subgroup is pseudo-Anosov.  The Farb--Mosher paper has inspired a continued interest in the convex cocompact subgroups of the mapping class group.  

Bowditch \cite{bowditch} showed that there are (at most) finitely many conjugacy classes of all pseudo-Anosov surface subgroups of $\Mod(S_g)$, for each $g$ (this result was also announced at about the same time by Groves and by Sapir).  So if such subgroups do exist, they are rare.

\bibliographystyle{plain}
\bibliography{problems}

\end{document}